\documentclass[10pt,a4paper]{article}
\makeatletter

\usepackage[left=3cm,top=3cm,bottom=3cm,right=3cm,head=1cm,foot=1cm]{geometry}
\usepackage[pdftex,bookmarks=true,bookmarksnumbered=true]{hyperref}
\usepackage{graphicx}
\usepackage{amssymb}
\usepackage{amsmath}
\usepackage{amsthm}
\usepackage{amsfonts}
\usepackage{amssymb}
\usepackage{amscd}
\usepackage{url}
\usepackage[all]{xy}

\newcommand{\ch}{\mathrm{ch}}

\newcommand{\sTr}{\mathrm{sTr}}
\newcommand{\h}{\textbf{H}}

\newcommand{\5}{\hspace{0,5cm}}

\newcommand{\3}{\vspace{0,3cm}}
\date{}
\title{Twisted K-theory constructions in the case of a decomposable Dixmier-Douady class}
\author{Antti J. Harju and Jouko Mickelsson\\
Department of Mathematics and Statistics, University of Helsinki}

\begin{document}

\maketitle

\begin{abstract}
	Twisted K-theory on a manifold $X$, with twisting in the 3rd integral cohomology,  is discussed in the case when 
$X$ is a product of a circle $\mathbb{T}$ and a manifold $M.$  The twist is assumed to be decomposable as a cup product of the basic integral one
form on $\mathbb{T}$ and an integral class in $H^2(M, \mathbb{Z}).$ This case was studied some time ago by V. Mathai, R. Melrose, and 
I.M. Singer. Our aim is to give an explicit construction for the twisted K-theory classes using a quantum field theory model, in the same spirit 
as the supersymmteric Wess-Zumino-Witten model is used for constructing (equivariant) twisted K-theory classes on compact Lie groups. \vspace{0,3cm}

\noindent Msc: 19L50, 53C08, 81T70 \vspace{0,3cm}

\end{abstract}

\section{Introduction} 

K-theory on a topological space $X$ can be twisted by an integral cohomology class $\sigma$ of degree 3. The class $\sigma$ can be either torsion
(this case was originally studied in \cite{DK}) or nontorsion \cite{Ro}.  The nontorsion case involves intrinsically infinite dimensional geometry since the
class $\sigma$ is the characteristic class of a principal bundle with the structure group $PU(\mathcal{H}),$ the projective unitary group of an infinite dimensional
separable complex Hilbert space $\mathcal{H}$. Partly because of this reason concrete constructions are available only in few cases. Best known of these is twisted K-theory
on a compact Lie group $G.$ It was shown by Freed, Hopkins, and Teleman \cite{FHT} that in the $G$ equivariant case the K-theory $K^*(G, \sigma)$ has a ring structure
isomorphic to the Verlinde ring in conformal field theory.  Concretely, the twisted -theory classes can be constructed from the quantized supersymmetric
Wess-Zumino-Witten model, \cite{Mi}. 

In this article we shall concentrate on the case $X= \mathbb{T} \times M,$ where $M$ is a compact manifold, $\mathbb{T}=\mathbb{T}_{\phi}=  S^1$  is a unit circle and the class $\sigma$ is represented as a product
\begin{eqnarray*}
\sigma = \frac{d\phi}{2\pi} \wedge \frac{\beta_M}{2 \pi i}
\end{eqnarray*}
 of the 1-form on $\mathbb{T}_{\phi}$ and a closed integral 2-form on $M.$ The Dixmier-Douady class associated to this type of product twistings was already studied by Raeburn and Rosenberg \cite{RR}, see Corollary 3.5 and Theorem 4.12,   and the twisted K-theory
 defined by the twisting  was  studied
by Mathai, Melrose, and Singer \cite{MMS}. (Actually, slightly more generally for the case of a product 
3-form, not necessarily arising from a factorization of the underlying space.) In particular, a formula for the Chern character  was derived in the decomposable case.  
The Chern character does not directly see the torsion classes in $K^*(X, \sigma).$  For this reason we want to analyze closer the torsion classes.
We also give a concrete formula for representatives of those classes using a quantum field theory construction similar to \cite{Mi} in the case of a compact
simply connected Lie group.  As a particular case, we have a construction for the (nonequivariant) torsion classes when $M$ is a torus. 

Our model is based on the methods in Hamiltonian quantization of fermionic fields. In general, consider a compact odd dimensional spin manifold $Z$. Suppose that $X$ is a compact manifold and that there is a family of Dirac operators, $\eth$, acting on the smooth spinors on $Z$ which is parametrized on $X$ and coupled to a connection of a vector bundle over $Z \times X$: for our purposes, it is sufficient to consider a trivial geometric fibration over $X$. Choose an open cover $\{U_a\}$ of $X$ together with real numbers $\{s_a\}$ so that $s_a \notin \text{spec}(\eth_A)$ for any $A \in U_a$. For each $x \in X$ we have a polarization of the Hilbert space of spinors onto the eigenspaces of $\eth_x$ with eigenvalues greated that $s_a$ and its orthogonal complement:
\begin{eqnarray*}
\mathcal{H}_x = \mathcal{H}_{x,a}^+ \oplus \mathcal{H}_{x,a}^-.
\end{eqnarray*}
Over each $U_a$ we can define a bundle of Fock space $\mathbf{F}_a$ which has the fibre at $x \in U_a$ given by 
\begin{eqnarray*}
\mathcal{F}_{x,a} = \bigwedge \mathcal{H}_{x,a}^+ \otimes (\bigwedge  \mathcal{H}_{x,a}^-)^*. 
\end{eqnarray*}
 We can fix the parametrization of $\{s_a\}$ so that $s_a \leq s_b$ if and only if $a < b$. Then locally, on the intersections $U_a \cap U_b = U_{ab}$, we have
\begin{eqnarray*}
 \textbf{F}_{a} \simeq \text{Det}_{ab} \otimes \textbf{F}_{b}   
\end{eqnarray*}
where $\text{Det}_{ab}$ is the line bundle whose fibre at $x \in U_{ab}$ is the determinant line of the finite dimensional $(s_a, s_b)$-spectral subspace of $\eth_x$. Then one can define a bundle of projective Fock space $\textbf{PF}$ over $X$ whose cohomology class is independent on the choice of the vacuum levels $\{s_a\}$. On the other hand, if the local bundles $\textbf{F}_a$ are given, then one would like to know if they can be glued together to a vector bundle over $X$. To this end we study triple overlaps of the cover and define local trivializations $f_{abc}: U_{abc} \rightarrow \mathbb{T}$ for the line bundles 
\begin{eqnarray*}
\text{Det}_{ab} \otimes \text{Det}_{bc} \otimes \text{Det}_{ac}^{-1}. 
\end{eqnarray*}
These trivializations define a cocycle in the {\v C}ech cohomology group  $[f] \in H^2(X, \underline{\mathbb{T}})$ which is the Dixmier-Douady cohomology class of the projective bundle $\textbf{PF}$. Indeed, the cocycle $[f]$ measures the obstruction to glue the local Fock bundles $\textbf{F}_a$ to a bundle over $X$. Recall that $H^2(X, \underline{\mathbb{T}})$ is isomorphic to $H^3(X, \mathbb{Z})$. Then the Dixmied-Douady class $[f]$ has a de Rham representative given by the degree 3-component of the characteristic class of the odd families index of $\eth$ over $X$, \cite{CMM}, \cite{Lott}.

We start the discussion in Section 2 by constructing a family of Dirac operators on a unit circle $Z = \mathbb{T}_{\theta}$ parametrized by $X= \mathbb{T}_{\phi} \times M$. 
As an input we need a complex line bundle $\lambda$ on $M$  with a connection $\nabla_M$ and curvature $\beta_M$. We shall introduce a Dirac family which has  spectral flow in the direction $\mathbb{T}_{\phi}$ and which is coupled to the line bundle $\lambda$ on $M$. The Atiyah-Singer families index theorem gives the Chern character
of this Fredholm family and one easily computes that the 3-form part gives the integral form $\frac{1}{4i\pi^2} d\phi \wedge \beta_M$. It then follows from the general principles in Hamiltonian quantization that there is a gerbe associated with this family. 

The construction of the gerbe is explained in Section 3 with details. This gerbe is realized as a projective Fock bundle over $\mathbb{T} \times M$. If pulled to the covering $\mathbb{R} \times M$ of $\mathbb{T} \times M$, the Dixmier-Douady class of the gerbe becomes trivial and there we get a true Fock bundle. In Section 4.1 we give a construction of a family of self-adjoint Fredholm operators parametrized by $\mathbb{R} \times M$ using a method similar to
\cite{Mi}. The supercharge family $Q_{(\phi, x)}$ (with $\phi\in\mathbb{R}$ and $x\in M$) transforms with respect to translations $\phi \mapsto \phi + 2\pi$ according to
$\hat{S}_x Q_{(\phi,x)} \hat{S}_x^{-1} =  Q_{(\phi- 2\pi,x)}$, for any given $x\in M$. However, there is no continuous map $x\mapsto \hat{S}_x$ from $M$ to unitary 
operators in the Fock space since the construction of $\hat{S}_x$ requires a choice of a unit vector in the fibre $\lambda_x$ of the complex line bundle $\lambda$. 
For this reason we obtain a genuinely twisted K-theory class on $X$ with twisting $\sigma  = \frac{1}{4i\pi^2} d\phi \wedge \beta_M$. We introduce another dependence on the topology of $M$ by tensoring the Fock vacua by a complex vector bundle $\xi$ over $M$. This has an effect on the twisted K-theory class of the family $Q$.

In Section 4.2 we compute the twisted K-theory groups subject to the decomposable twist $\sigma$. The group $K^1(\mathbb{T} \times M, \sigma)$ is isormorphic to a group extension of  
\begin{eqnarray}\label{K1}
	\{x \in K^1(M): x = x \otimes \lambda\} \5 \text{by} \5 \frac{K^{0}(M)}{K^{0}(M) \otimes (1 - \lambda)}.
\end{eqnarray}
Torsion subgroups are common in twisted K-theory: for example in the case of a unit sphere  $M= S^2$ we get the torsion elements in $K^1(X, \sigma) = \mathbb{Z} \oplus \mathbb{Z}_k$ where $k$ is the degree of $\lambda$. In Section 5, differential twisted K-theory is applied to study how our approach for the analytic realization of the twisted K-theory corresponds to the torsion subgroups. Finally, in Section 6 we compute a Chern character for the index using the Quillen superconnection formalism. The Chern character is well defined as a differential form on the covering space $\mathbb{R} \times M$. One can try to pull this character to a differential form on $\mathbb{T} \times M$. However, now the Chern character is defined only modulo multiplication by $\text{ch}(\lambda)$, the Chern character of the twisting line bundle. We shall intoduce a new target group which is essentially the odd de Rham group $H^{\text{odd}}(\mathbb{T} \times M)$ quotioned by the normal subgroup whi!
 ch kills the normal subgroup $K^0(M) \otimes (1 - \lambda)$ in \eqref{K1} on the level of the de Rham cohomology. This construction can be applied to extract information of the analytic realizations of $K^1(\mathbb{T} \times M, \sigma)$ including torsion components. An explicit representative for the character of the supercharge $Q$ is derived. \vspace{0.3cm} 

\noindent \textbf{Conventions.} The symbol $\mathcal{H}$ denotes a separable complex Hilbert space. The unitary group $U(\mathcal{H})$ with the norm topology is contractible by Kuiper's theorem \cite{K}. This topology is not suitable for our constructions since we are employing representations of loop groups, which are not continuous in this topology. However, $U(\mathcal{H})$ is also contractible in the strong operator topology \cite{DD} and this is the topology we are going to use. $PU(\mathcal{H})$ denote the group of projective unitary transformations. 

We follow Atiyah and Segal in the definition of twisted K-theory on a compact space $X$, \cite{AtSe}.  The ordinary K-theoretic classification spaces are not well behaving in the twisted setup since we need to have a continuous action of $U(\mathcal{H})$ on them. Consider the space of bounded operators $\mathcal{B}(\mathcal{H})$ equipped with the strong operator topology and the space of compact operaotors $\mathcal{K}(\mathcal{H})$ with the norm topology. Let $\textbf{Fred}^{(0)}$ denote the subspace in $\mathcal{B}(\mathcal{H})$ of the operators with a parametrix: $A \in \textbf{Fred}^{(0)}$ if there is an operator $B \in \mathcal{B}(\mathcal{H})$, such that $1 - AB$ and $1 - BA$ are compact. Let $\textbf{Fred}^{(1)}$ be the subspace of self-adjoint operators in $\textbf{Fred}^{(0)}$ which have both positive and negative essential spectrum. Let $X$ be a compact space and $\sigma \in H^3(X, \mathbb{Z})$. Suppose that, under the standard isomoprhism 
\begin{eqnarray*}
H^1(X, \underline{PU}(\mathcal{H})) \rightarrow H^3(X, \mathbb{Z}),
\end{eqnarray*}
the class of $\sigma$ is realized by a $PU(\mathcal{H})$-bundle $\textbf{P}$. Associated to the principal bundle $\textbf{P}$ there is a projective bundle of Hilbert spaces. Under the adjoint action of $PU(\mathcal{H})$ on $\textbf{Fred}^{*}$, there are the bundles of Fredholm operators for $* = 0,1$. The twisted K-theory groups are defined to be the groups of homotopy classes of continuous sections in these bundles
\begin{eqnarray*}
K^{*}(X, \sigma) = \{[F]: F \in \Gamma(\textbf{P} \times_{\text{Ad}} \textbf{Fred}^{*})\}.
\end{eqnarray*}
The above definition of twisted K-theory is used in the local formulation: choose a trivializing open cover $\{U_i\}$ of $X$ for the bundle $\textbf{P}$. Let $g_{ij}: U_i \cap U_j \to PU(H)$ be the corresponding transition functions. Then a twisted K-theory element is given in terms of local functions $F_i : U_i \to \textbf{Fred}^{*}$ such that $F_i = g_{ij} F_j g_{ij}^{-1}$ on the overlaps. 

Our analysis is based on unbounded Fredholm operators. The spaces $\textbf{Fred}^*_{\Psi}$ with $* = 0,1$ are equipped with the topology which is induced by the map 
\begin{eqnarray*}
\Psi (Q) = \frac{Q}{\sqrt{1 + Q^2}} 
\end{eqnarray*}
and the elements in $\textbf{Fred}^0_{\Psi}$ are the unbounded Fredholm operators which map to an element in $\textbf{Fred}^{(0)}$ under $\Psi$, and similarly for $\textbf{Fred}^1_{\Psi}$.

\section{The Dixmier-Douady class and the families index theorem} 

Consider a compact product manifold $\mathbb{T}_{\phi} \times M$ with a nontrivial decomposable integral $3$-cohomology class, \cite{MMS}, represented in the de Rham cohomology by
\begin{eqnarray*}
\sigma = \frac{d \phi}{2 \pi} \wedge \frac{\beta_M}{2 \pi i},
\end{eqnarray*}
where $\beta_M$ represents a class of $H^2(M, 2 \pi i \mathbb{Z})$. We are interested in $K$-theory twisted by a gerbe and therefore we can exploit the Hamiltonian quantization to build a gerbe over $\mathbb{T}_{\phi} \times M$. The first goal is to construct a family of Dirac operators on $\mathbb{T}_{\phi} \times M$ with a three form component in its index given by the decomposable cocycle $\sigma$. \vspace{0.3cm}

\noindent \textbf{2.1.} Consider a $2$-torus $\mathbb{T}^2$ with angle variables $(\theta, \phi)$. We choose an open cover $\{\mathbb{T}_+, \mathbb{T}_- \}$ for $\mathbb{T}$ such that $\mathbb{T}_{+-} = \mathbb{T}_+ \cap \mathbb{T}_-$ consists of two disconnected arcs, one which is a neighbourhood of $-1$ and another a neighbourhood of $1$. We denote these by $\mathbb{T}_{+-}^{(-1)}$ and $\mathbb{T}_{+-}^{(1)}$. 

The isomorphism classes of line bundles over $\mathbb{T}^2$ are classified by $\mathbb{Z}$ since $H^2(\mathbb{T}^2, \mathbb{Z}) = \mathbb{Z}$.  The bundle
$\lambda_1$ corresponding to a generator of the cohomology group can be described as follows: if $\psi$ is a smooth section of $\lambda_1$, then $\psi(\theta, \phi + 2\pi)
= e^{ i \theta} \psi(\theta,\phi)$. After puling back with the map $\mathbb{R} \times \Bbb T \to \Bbb T \times \Bbb T,$ sending  $\phi$ to $\phi$ mod $2\pi,$ a connection of this bundle can be defined by
\begin{eqnarray*}
\nabla_1 =  d\theta \otimes \partial_{\theta}+ d\phi\otimes \partial_{\phi} -  \frac{i}{2\pi} d\theta \otimes \phi.
\end{eqnarray*}
The curvature of the connection is the cocycle in de Rham cohomology
\begin{eqnarray*}
\nabla_1^2 = \frac{ i}{2\pi}  d\theta \wedge d\phi \in H^2(\mathbb{T}^2,2\pi i \Bbb Z).
\end{eqnarray*}

Consider a smooth manifold $M$ with nontrivial second cohomology and fix a line bundle $\lambda$ with a connection so that the curvature is equal to $\beta_M$ whose class in $H^2(M)$ is nontrivial. Now $\tilde \lambda = \lambda_1 \boxtimes \lambda$ defines a line bundle over $\mathbb{T}^2 \times M$.
Consider a smooth fibration
\begin{eqnarray*}
\mathbb{T}_{\theta} \hookrightarrow \mathbb{T}_{\theta} \times
\mathbb{T}_{\phi} \times M \rightarrow \mathbb{T}_{\phi} \times M.
\end{eqnarray*}
At each $(\phi, x) \in \mathbb{T}_{\phi} \times M$, the bundle $\tilde \lambda$ restricted to the fibre $\mathbb{T}_{\theta}$ defines a line bundle $\lambda(\phi, x) \rightarrow \mathbb{T}_{\theta}$. In fact, the sections of this bundle are periodic in the direction $\theta$ and therefore the bundle $\lambda(\phi, x)$ is the product $\mathbb{T}_{\theta} \times \mathbb{C}$.

At each point $(\phi, x)$ we define a Hilbert space $\mathcal{H}(\phi, x) = L^2(\mathbb{T}_{\theta}, \lambda(\phi,x))$ of $L^2$-sections on $\mathbb{T}_{\theta}$ with values in $\lambda(\phi,x)$. Then
\begin{eqnarray*}
\textbf{H} = \coprod_{(\phi,x) \in \mathbb{T} \times M}\mathcal{H}(\phi,x)
\end{eqnarray*}
is a locally trivial bundle of Hilbert spaces over $\mathbb{T} \times M$. As a Hilbert bundle it is trivial by Kuiper's theorem. However, considered as a reduced bundle with the structure group of smooth $\mathbb{T}$ valued gauge transformations, the group $L\mathbb{T}_{\theta}$ of smooth endomorphism of $\mathbb{T}_{\theta},$ it is nontrivial.
The gauge group acts on each fibre $\mathcal{H}(\phi,x)$ by multiplication: $m: L\mathbb{T}_{\theta} \times  \mathcal{H}(\phi,x) \rightarrow \mathcal{H}(\phi,x)$. The group $\mathbb{Z}$ of translations over $\mathbb{T}_{\phi}$ acts on the sections of $\h$ from the left by the rule 
\begin{eqnarray*}
a.\varphi(\phi,x) = m(e^{ i a \theta}) \varphi(\phi,x),
\end{eqnarray*}
for all $a \in \mathbb{Z}$ and $(\phi, x) \in \mathbb{T}_{\phi} \times M$. 

Associated to the isomorphism class of the line bundle $\lambda$, there is a {\v C}ech cohomology class in $H^1(M, \underline{\mathbb{T}})$. Fix a good open cover $\{V_i, 1 \leq i \leq n \}$ for $M$. Then we can fix a representative $h$ of the class of $\lambda$ in $H^1(M, \underline{\mathbb{T}})$ which has the components $h_{ij} \in C^{\infty}(V_{ij}, \Bbb T)$ where $V_{ij} = V_i \cap V_j$. We fix a cover $\{U_a: 1 \leq a \leq 2n\}$ for $\mathbb{T}_{\phi}\times M$
such that
\begin{eqnarray*}
U_{i} = \mathbb{T}_+ \times V_i, \5 U_{i + n } = \mathbb{T}_- \times V_i, \5 1 \leq i \leq n. 
\end{eqnarray*}
The transition functions of $\h$ define a totally antisymmetric cocycle $g$ whose cohomology class lies in $H^1(\mathbb{T} \times M, \underline{L\mathbb{T}_{\theta}})$. The components of this cocycle are the locally defined functions $g_{ab}:  U_{ab} \rightarrow L\mathbb{T}_{\theta}$ for $1 \leq a,b  \leq  2n$ and the nonidentity components are determined by the antisymmetry in the indices $a,b$ and the relations
\begin{eqnarray}\label{Hdata}
 g_{ij}(\phi,x) &=& h_{ij}(x),\nonumber \\  
 g_{i+n, j+n}(\phi,x) &=& h_{ij}(x),\nonumber  \\
 g^{(-1)}_{i,j+n}(\phi,x) &=& h_{ij}(x),\nonumber  \\ 
 g^{(1)}_{i,j+n}(\phi,x) &=& m(e^{ i  \theta}) h_{ij}(x)
\end{eqnarray}
for $1 \leq i,j \leq n$ and $g^{(\pm 1)}_{i,j+n}$ is defined in the open set $\mathbb{T}^{(\pm 1)}_{+-} \times V_{ij}$.
These transition functions satisfy the $1$-cocycle condition $g_{ab} g_{bc} g_{ac}^{-1} = 1$ in their domains. Then we can write
\begin{eqnarray*}
\h = \Big[\coprod_a U_a \times \mathcal{H}\Big] / \sim
\end{eqnarray*}
where we use the equivalence relation in the overlaps: $(\phi,v) \sim (\phi', w)$ if and only if $\phi = \phi'$ in $U_{ab}$ and $w = g_{ab} v$ in
$\mathcal{H}$. \3

\noindent \textbf{2.2.} The free Dirac operator $-i \partial_{\theta}$ is an unbounded self-adjoint operator on each fibre
$\mathcal{H}(\phi,x)$. The space of vector potentials on each fibre is given by
$\mathcal{A} \simeq C^{\infty}(\mathbb{T}_{\theta}) \otimes i\mathbb{R}$. The gauge group $L
\mathbb{T}_{\theta}$ acts on the Dirac operators by conjugation, leading to the action $A\mapsto A + g^{-1} dg$ on gauge potentials. The gauge orbit space is
$\mathcal{A} / L \mathbb{T}_{\theta}$ which can be identified with a
circle. Thus, $\mathbb{T}_{\phi}$ has a natural interpretation of a space of
gauge potentials which we twist with the bundle $\lambda$ on $M$. Actually, it is sufficient to consider constant vector potentials $\phi$
parametrized by the real line $\Bbb R.$ The gauge transformations by $\mathbb{T}$ valued functions $e^{i\theta}$ on $\mathbb{T}_{\theta}$ change the parameter
$\phi \mapsto \phi +2\pi,$ so again the family $-i\partial_{\theta} +\frac{\phi}{2\pi} $ modulo gauge transformations is parametrized by $\mathbb{R}/2\pi \mathbb{Z} = \mathbb{T}_{\phi}.$
After twisting this family by the line bundle $\lambda$ over $M$ we get a family parametrized by $X= \mathbb{T}_{\phi} \times M.$  The Dirac family
\begin{eqnarray*}
\eth: (\phi, x) \mapsto - i \partial_{\theta} + \frac{\phi}{2 \pi}
\end{eqnarray*}
 is twisted by the complex line bundle over $\Bbb T^2 \times M$ with connection $\nabla_1 \otimes \nabla_M$ and the total curvature
\begin{eqnarray*}
F =  \frac{i}{2\pi} d\theta \wedge d\phi
 +  \beta_M \in H^2(\mathbb{T}_{\theta}\times X, 2\pi i \Bbb Z).
\end{eqnarray*}

The Dirac family $\eth$ defines an eigenvalue problem at each $(\phi,x) \in \mathbb{T}_{\phi} \times M$. If we let the angle $\phi$ vary from
$0$ to $2\pi$, then there is a translation in the set of eigenvalues as they all
grow by $2 \pi$. Because of the spectral flow, the group element of 
$K^1(\mathbb{T}_{\phi} \times M)$ defined by the Fredholm family is nontrivial. In fact, the
spectral flow produces a nontrivial cocycle of $H^1(\mathbb{T}_{\phi} \times M,
\mathbb{Z})$ via the index map. The twisting bundle $\lambda$
produces another nontrivial class, a three form in
$H^{\mathrm{odd}}(\mathbb{T}_{\phi} \times M, \mathbb{Z})$. The local
index formula, $\mathrm{ind}: K^1(\mathbb{T}_{\phi} \times M) \rightarrow
H^{\mathrm{odd}}(\mathbb{T}_{\phi} \times M)$, gives 
\begin{eqnarray*}
\mathrm{ind}(\eth) &=& \int_{\mathbb{T}_{\theta}}  \ch(\lambda_1 \boxtimes
\lambda) \\ &=& \int_{\mathbb{T}_{\theta}}  \exp(\frac{\nabla_1^2}{2 \pi i}) \wedge
\exp(\frac{\beta_M}{2 \pi i}) \\
&=& \int_{\mathbb{T}_{\theta}} \exp\left(  \frac{1}{4\pi^2} d\theta\wedge d\phi   +  
 \frac{ \beta_M }{2 \pi i} \right) \\
&=& \frac{d\phi}{2\pi}   + \frac{d\phi}{2\pi}   \wedge \frac{ \beta_M }{2 \pi i} + \cdots .
\end{eqnarray*}
The A-roof genus  on $\mathbb{T}^2\times M$ does not contribute on this level in the character formula.
The three cohomology part is exactly the decomposable $3$-cohomology class.

\section{Projective Fock bundle} 

\noindent \textbf{3.1.} Let $\mathcal{H}$ be a complex separable Hilbert space. The algebra $A$ is called a canonical anticommutation relations (CAR) algebra over $\mathcal{H}$ if
there is an antilinear mapping $a: \mathcal{H} \rightarrow A$ such that ${a(u): u \in \mathcal{H}}$ generate a unital $C^*$-algebra $A$ which fulfills
\begin{eqnarray*}
[a(u), a(v)]_+ = 0 \5 \hbox{and}\5 [a(u), a(v)^* ]_+ = \langle u, v \rangle 1
\end{eqnarray*}
where $[\cdot,\cdot]_+$ denotes the anticommutator, for all $u,v \in \mathcal{H}$. The CAR algebra is unique up to $C^*$-algebra isomorphism. 

For a fixed  $(\phi,x) \in \mathbb{T} \times M$, the Dirac operator $\eth_{\phi,x}$ defines a polarization $\mathcal{H} = \mathcal{H}^+ \oplus \mathcal{H}^-$ such that $\mathcal{H}^+$ is spanned by the nonnegative eigenstates. A Fock space $\mathcal{F}$ is a Hilbert space with a
vacuum vector $| 0 \rangle$ and the CAR algebra acts on the vacuum such that
\begin{eqnarray*}
a(u)|0 \rangle = 0 = a^*(v)|0 \rangle\5 \hbox{for all}\5  u \in \mathcal{H}^+, v \in \mathcal{H}^-, 
\end{eqnarray*}
and the basis of a Fock space is spanned by 
\begin{eqnarray*}
a(u_{i_1}) \cdots a(u_{i_k})
a^*(u_{j_1}) \cdots a^*(u_{j_l}) |0 \rangle, \5 \hbox{for}\5 u_{i_{\nu}} \in
\mathcal{H}^-, u_{j_{\nu}} \in \mathcal{H}^+.
\end{eqnarray*}
We can think of the vacuum as the formal infinite wedge product
\begin{eqnarray*}
|0 \rangle =  u_{-1} \wedge u_{-2} \wedge u_{-3} \wedge \cdots
\end{eqnarray*}
and the general basis vector as
\begin{eqnarray*}
u_{j_1} \wedge u_{j_2} \wedge u_{j_3} \wedge \cdots
\end{eqnarray*}
where $j_1 > j_2 > j_3> \cdots$ are integers such that all the negative integers  except a finite number are included in the sequence. The representation of CAR is irreducible. There exists a densely defined charge operator $N$ which acts on a basis vectors by
\begin{eqnarray*}
N a(u_{i_1}) \cdots a(u_{i_k}) a^*(u_{j_1}) \cdots a^*(u_{j_l}) |0 \rangle = (l-k) a(u_{i_1}) \cdots a(u_{i_k}) a^*(u_{j_1}) \cdots a^*(u_{j_l}) |0 \rangle
\end{eqnarray*}
and its action is extended on $\mathcal{F}$ by linearity. The Fock space can be presented as a completion of the algebraic direct sum over the charge subspaces
\begin{eqnarray*}
\mathcal{F} = \hat{\bigoplus_{k \in \mathbb{Z}}} \mathcal{F}^{(k)}
\end{eqnarray*}
where $\mathcal{F}^{(k)}$ is the subspace of charge $k$: the eigenspace of $N$ for the eigenvalue $k$. 

In the group $L \mathbb{T}_{\theta}$  of smooth loops in $\mathbb{T}_{\theta}$ any element is of the form $e^{2 \pi iF}$ such that $F: \mathbb{R} \rightarrow
\mathbb{R}$ is a smooth function and 
\begin{eqnarray*}
F(\theta + 2\pi) = F(\theta) +  n_F.
\end{eqnarray*}
$n_F \in \mathbb{Z}$ is the winding number of the loop. Then $f(\theta) = F(\theta) - \frac{n_F \theta}{2\pi}$ is invariant under the translations $\theta\mapsto \theta +2\pi$ and thus it can be
expanded as a Fourier series $f = \sum f_k e_k$, where $f_k$ are the Fourier coefficients for all $k \in \mathbb{Z}$. Since $f$ is real valued these satisfy $\overline{f_{k}} = f_{-k}$. We can write 
\begin{eqnarray*}
L \mathbb{T}_{\theta} = SL \mathbb{T}_{\theta} \times C \mathbb{T}_{\theta}
\end{eqnarray*}
 such that the charge subroup $C \mathbb{T}_{\theta}$ consists of the group elements $e^{ 2 \pi i f_0 + i n_F
\theta}$ and $SL \mathbb{T}_{\theta}$ consists of $e^{2 \pi i\sum_{k \neq 0} f_k e_k }$.

The loop group $L \mathbb{T}_{\theta}$ is a subgroup of the restricted unitary group $U_{\mathrm{res}}(\mathcal{H}^+ \oplus \mathcal{H}^-)$ which has a
positive energy representation on a Fock space, \cite{PS}. The action of $U_{\mathrm{res}}(\mathcal{H}^+ \oplus \mathcal{H}^-)$ can be implemented on the
Fock space as a projective representation such that
\begin{eqnarray*}
U(g) a(u) U(g^{-1}) = a(g.u), \5 U(g) a^*(v) U(g^{-1}) = a^*(g.v)
\end{eqnarray*}
for all $g \in U_{\mathrm{res}}(\mathcal{H}^+ \oplus \mathcal{H}^-)$ and $u,v \in \mathcal{H}$. The subroup $SL \mathbb{T}_{\theta}$ lies in the
connected component of the identity of $U_{\mathrm{res}}(\mathcal{H}^+ \oplus \mathcal{H}^-)$ and each charge subspace is invariant under this action.  The
subgroup $C \mathbb{T}_{\theta}$ has infinitely many disconnected components labeled by $n_F \in \mathbb{Z}$. If $e^{2 \pi i f_0 + i n_F \theta} \in C \mathbb{T}_{\theta}$,
then it is represented on the Fock space by
\begin{eqnarray}\label{gerbe1}
U(e^{2 \pi i f_0 + i n_F \theta}) = e^{\pi i f_0 N} S^{n_F} e^{\pi i f_0 N},
\end{eqnarray}
where $S$ is a shift operators which sends each charge subspace $\mathcal{F}^{(k)}$ to $\mathcal{F}^{(k+1)},$ that is, 
\begin{eqnarray*}
SNS^{-1} = N-1.
\end{eqnarray*}
More precisely, given an orthonormal basis $\{u_n : n \in \mathbb{Z} \}$ for $\mathcal{H}$, then 
\begin{eqnarray*}
S|0 \rangle &=& a^*(u_0) |0 \rangle, \\
S^{-1}|0 \rangle &=& a(u_{-1}) |0 \rangle, \\
S a(u_n) S^{-1} &=& a(u_{n+1}), \\ S a^*(u_n) S^{-1} &=& a^*(u_{n+1}). 
\end{eqnarray*}
The positive energy representation of $L \mathbb{T}_{\theta}$ are projective: there is a group $2$-cocycle $c: L \mathbb{T}_{\theta} \times L \mathbb{T}_{\theta} \rightarrow
\mathbb{T}$ such that the unitary representation satisfies 
\begin{eqnarray}\label{proj}
U(e^{iF}) U(e^{iG})  = U(e^{i(F+G)})c(e^{iF}, e^{iG}) .
\end{eqnarray}
We denote by $\mathcal{PF}$ the projective Fock space $\mathcal{F}/\mathbb{T}$. Then $U$ defines a representation of $L\mathbb{T}_{\theta}$ on $\mathcal{PF}$.\3

\noindent \textbf{3.2.} Next we construct a projective Fock bundle associated to the Dirac family in 2.2.  Fix a complex line bundle $\lambda$ over $M$ and a cover $\{V_i\}$ of $M$ which trivializes $\lambda$, as in 2.1. Then $\lambda$ is extended on the space $\mathbb{T} \times M$ so that the new transition functions satisfy $h_{ab}(\phi,x) = h_{ab}(x)$ for all $(\phi,x) \in \mathbb{T} \times M$.  Similarly we can extend $\lambda$ to $\mathbb{R} \times M$. On $\mathbb{T} \times M$ we use the open cover $\{U_a: 1 \leq a \leq 2n\}$, recall 2.1. We also write
\begin{eqnarray*}
U_{i,j+n}^{(1)} = \mathbb{T}^{(1)}_{+-} \times V_{ij}
\end{eqnarray*}
where $\mathbb{T}^{(1)}_{+-}$ is the neighborhood of $1$ in the unit circle. The nonidentity components of the gerbe cocycle will localize on these components. 

The spectrum of the Dirac operator $\eth_{(\phi,x)}$ is the set $\{n + \frac{\phi}{2 \pi}: n \in \mathbb{Z}\}$. The corresponding eigenstates are $u_n := \frac{1}{\sqrt{2\pi}} e^{i n \theta}$. Then one can easily fix vacuum parameters $s_+$ and $s_-$ so that locally, over $\mathbb{T}_+ \times M$ and $\mathbb{T}_- \times M$, the vacuums of the local Fock bundles $\textbf{F}_{\pm}$ are
\begin{eqnarray*}
|0 \rangle_{\pm} = u_{-1} \wedge u_{-2} \wedge u_{-3} \wedge \cdots
\end{eqnarray*}
Recall, \eqref{Hdata}, that over the intersection $\mathbb{T}_{+-}^{(1)} \times M$, the Hilbert space state $u_k$ at $(\phi + 2 \pi, x) \in \mathbb{T}_- \times M$ is identified with $u_{k+1}$ at $(\phi, x) \in \mathbb{T}_+ \times M$. The corresponding transformation implemented to the Fock spaces creates a state from the vacuum and in this process the charge raises by one. The Hilbert bundle $\textbf{H}$ over $\mathbb{T} \times M$ has been tensored by the line bundle $\lambda$. Therefore, the charge $k$ sector of the Fock bundle behaves on the submanifold $M$ as the line bundle $\lambda^{\otimes k}$. The process to create a Fock state is thus defined only up to a phase: one needs to fix a phase on each fibre of $\lambda$ over $M$. Therefore we start by considering a projectivization of the Fock bundle. The projective bundle can be described by starting with the trivial bundle with projective Fock space fibres $\mathbf{PF}_0$ over the covering space $\mathbb{R} \times M$. The!
 n we define
\begin{eqnarray*}
	\textbf{PF} = \textbf{PF}_0 / \sim 
\end{eqnarray*}
where $\sim$ is the equivalence relation $(\phi,x, \Psi) \sim (\phi',x',\Psi')$ if and only if $\phi' = \phi + 2\pi n$ in $\mathbb{R}$, $x=x'$ in $M$ and $\Psi' = S_{x}^n \Psi$ in the fibre for all $n \in \mathbb{Z}$. We have applied to projective operator family $S: M \rightarrow PU(\mathcal{H})$. Now $\textbf{PF}$ is a projective bundle on $\mathbb{T} \times M$ and its Dixmier-Douady class is determined by a lift of the transition functions $S$ to unitary operators. 

Next we proceed to give an explicit description for the Dixmier-Douady class in terms locally defined data. The process to raise the charge of the Fock states under translations will be localized on the intersection $\mathbb{T}_{+-}^{(1)} \times M$. So on this intersection there is the bundle isomorphism 
\begin{eqnarray*}
\textbf{F}_+ \simeq \lambda \otimes \textbf{F}_-.
\end{eqnarray*}
However, there is no way to define a unitary family over $\mathbb{T}_{+-}^{(1)} \times M$ that would create a Fock state of topological type $\lambda$ simply by lifting the projectively defined shift operators $S$ to a unitary family on $M$: One needs to fix phases for $\lambda$ in every fibre but $\lambda$ is nontrivial so this cannot be done continuously over $M$. Instead, one can fix local frames $l_i: V_i \rightarrow \mathbb{T}$ and then lift the shift operators over $U^{(1)}_{i,j+n}$ with $1 \leq i,j \leq n$ accordingly. This results the local bundle  isomoprhisms
\begin{eqnarray*}
	|0 \rangle_{i} \simeq  \lambda_{i,j+n} \otimes |0 \rangle_{j+n}
\end{eqnarray*}
where $|0 \rangle_a$ is the vacuum subbundle on $U_a$ for all $1 \leq a \leq 2n$, and $\lambda_{i,j+n}$ is the bundle $\lambda$ on $U^{(1)}_{i,j+n}$ trivialized over $U_i$. In the inverse relations we have $\lambda_{j+n,i} = \lambda_{i,j+n}^{-1}$ which is the bundle $\lambda^{-1}$ on $U^{(1)}_{i,j+n}$ trivialized over $U_i$. The Dixmier-Douady class of $\textbf{PF}$ is represented by the cocycle $f$ which defines a class in $H^2(\mathbb{T} \times M, \underline{\mathbb{T}})$ and whose components $f_{abc}$ are totally antisymmetric and determined by a choice of a trivialization of  
\begin{eqnarray*}
	\lambda_{ab} \otimes \lambda_{bc} \otimes \lambda_{ac}^{-1} \5 \text{for all}\5 1 \leq a,b,c \leq 2n.
\end{eqnarray*}
In our case, up to permutation of the indices, only the local bundles of type 
\begin{eqnarray}\label{vacuums}
	\lambda_{ij} \otimes \lambda_{j,k+n} \otimes \lambda_{i,k+n}^{-1} \5 \text{and} \5 \lambda_{i,j+n} \otimes \lambda_{j+n,k+n} \otimes \lambda_{i,k+n}^{-1}
\end{eqnarray}
for $1 \leq i,j,k \leq n$ can contribute to the Dixmier-Douady class because there are no change of the vacuum type associated to other combinations. In fact the bundles on the right side of \eqref{vacuums} is canonically trivial over $U_{i}$. For the bundle on the left side, we can use the transition functions $h_{ij}$ to transform the component $\lambda_{j,k+n}$ to the trivialization over $U_i$. The component $\lambda_{ij}$ is a canonically trivial bundle whereas $\lambda_{i,k+n}^{-1}$ is trivialized over $U_i$. Thus, we have 
\begin{eqnarray*}
f^{(1)}_{i,j,k+n} : U^{(1)}_{i,j,k+n} \rightarrow \mathbb{T}; \5 f^{(1)}_{i,j,k+n}(\phi,x) = h_{ij}(x).
\end{eqnarray*}
Similarly one checks that these components are antisymmetric under the permutations of the indices $i,j,k+n$. So far we have proved the following. \3

\noindent \textbf{Proposition.} The Dixmier-Douady class of the projective bundle $\textbf{PF}$ defines a class in the {\v C}ech cohomology group $H^2(\mathbb{T} \times M, \underline{\mathbb{T}})$ which is represented by the cocycle $f$ whose nonidentity components are determined by the total antisymmetry under permutations of its indices and by the relations 
\begin{eqnarray}\label{f}
f^{(1)}_{i,j,k+n}(\phi,x) = h_{ij}(x), \5 1 \leq i ,j,k \leq n.
\end{eqnarray}
The representative of $f$ in $H^3(\mathbb{T} \times M, \mathbb{R})$ lies in the cohomology class of the decomposable cocycle
\begin{eqnarray*}
\sigma = \frac{d\phi}{2\pi} \wedge \frac{\beta_M}{2 \pi i} = \text{ind}(\eth)_{[3]}
\end{eqnarray*}
in the de Rham cohomology complex. \vspace{0,3cm}

\noindent \textbf{3.3.} Next we write down transition functions for a gerbe with Dixmier-Douady class equivalent to $f$ in the {\v C}ech cohomology group $H^2(\mathbb{T} \times M, \underline{\mathbb{T}})$. We apply the representation theory of $L \mathbb{T}$ on the projective Fock space and then proceed with the usual bundle reconstruction. In 2.1 we defined the transition data for the Hilbert bundle $\textbf{H}$, recall \eqref{Hdata}. We let these cocycles act on the projective Fock spaces $\mathcal{PF}$ under the projective representation $U$ defined in \eqref{gerbe1}. Now $S$ is defined up to a phase since we need to choose a basis for the Fock states it creates. Then we define a projective Fock bundle 
\begin{eqnarray*}
	\textbf{PF}' = \Big[ \coprod_{a} U_a \times \mathcal{PF} \Big] / \sim
\end{eqnarray*}
where the equivalence relation is associated to the cocycle $U(g_{ab}): U_{ab} \rightarrow PU(\mathcal{F})$. The lifted transition functions $\hat{g}_{ab}: U_{ab} \rightarrow U(\mathcal{F})$ are determined by the antisymmetry of the indices and  
\begin{eqnarray}\label{Fdata}
 \hat{g}_{ij}(\phi,x) &=& (h_{ij}(x))^N, \nonumber \\  
 \hat{g}_{i+n, j+n}(\phi,x) &=& (h_{ij}(x))^N, \nonumber \\
 \hat{g}^{(-1)}_{i,j+n}(\phi,x) &=& (h_{ij}(x))^N, \nonumber \\ 
 \hat{g}^{(1)}_{i,j+n}(\phi,x) &=& (h_{ij}(x))^{\frac{N}{2}} S (h_{ij}(x))^{\frac{N}{2}}
\end{eqnarray}
for all $1 \leq i,j \leq n$. The basis of the Fock space $\mathcal{F}$ is fixed and so $S$ is defined as a bounded operator on $\mathcal{F}$ but in this formulation $S$ does not depend on the point $x \in M$. The depdendence on the base is now associated to the transition functions $(h_{ij}(x))^N$ which operate on a charge $k$ state over $x$ exactly as the transition functions of $\lambda^{\otimes k}$ do. Then there is a cocycle $f'$ representing a class of $H^2(\mathbb{T} \times M, \underline{\mathbb{T}})$ defined by $f'_{abc} = \hat{g}_{ab} \hat{g}_{bc} \hat{g}_{ac}^{-1}$. We apply the relations
\begin{eqnarray}\label{gerbe2}
S (h_{ab}(\phi,x))^{\frac{N}{2}} S^{-1} = (h_{ab}(\phi,x))^{\frac{N-1}{2}}
\end{eqnarray}
in $U_{ab}$ for all $1 \leq a,b \leq 2n$. Applying \eqref{Fdata}, \eqref{gerbe2} and the relations $(h_{ij} h_{jk} h_{ki})^N = 1$ for all $i,j,k$ we find a totally antisymmetric {\v C}ech cocycle  with the nonidentity components determined by
\begin{eqnarray*}
(f')^{(1)}_{i,j,k+n}(\phi,x) = h_{ij}^{1/2}(x) \5 \text{and} \5 (f')^{(1)}_{i,j+n,k+n}(\phi,x) =  h_{jk}^{-1/2}(x). 
\end{eqnarray*}
for $1 \leq i,j,k \leq n$. \vspace{0,3cm}

\noindent \textbf{Proposition.} The Dixmier-Douady class of $\textbf{PF}'$ is represented by $f'$ and in the {\v C}ech cohomology group $H^2(\mathbb{T} \times M, \underline{\mathbb{T}})$ it is equivalent to the cocycle $f$ defined in \eqref{f} and to the totally antisymmetric cocycle $f''$ with the nontrivial components determined by 
\begin{eqnarray}\label{gerbecocycle2}
(f'')^{(1)}_{i,j+n,k+n}(\phi,x) = (h_{jk})^{-1}(x), \5 1 \leq i,j,k \leq n.
\end{eqnarray}

\noindent Proof. We can use the property \eqref{proj} to adjust the shift operators $S$ on the right side of the transition functions $(h_{ij})^N$ in the definition of each $\widehat{g}_{i,j+n}$ for $1 \leq i,j \leq n$. Then we can write
\begin{eqnarray*}
(f')^{(1)}_{i,j+n,k+n} &=& (h_{ij})^N S (h_{jk})^N S^{-1} (h_{ik})^{-N} c(h_{ij}, e^{i \theta})^{-1} c(h_{ik}, e^{i \theta}) \\
&=& h_{jk}^{-1} c(h_{ij}, e^{i \theta})^{-1} c(h_{ik}, e^{i \theta}).
\end{eqnarray*}
The $\mathbb{T}$ valued functions $c(h_{ij}, e^{i \theta})^{-1} c(h_{ik}, e^{i \theta})= (\delta c^{-1})_{i,j+n,k+n}$ are components of the coboundary $\delta c^{-1}$ and $c^{-1}$ is the cochain with the nonidentity components 
\begin{eqnarray*}
	c^{-1}_{i,j+n}(\phi,x) = c(h_{ij}(x), e^{i \theta})^{-1} \5 \text{and} \5 c^{-1}_{j+n,i}(\phi,x) = c(h_{ij}(x), e^{i \theta})
\end{eqnarray*}
in $U^{(1)}_{i,j+n}$  for all $1 \leq i,j \leq n$. Similarly one can check that the permutations of these components are equal to the components of 
\eqref{gerbecocycle2} times the coboundary of $c^{-1}$. Moreover, 
\begin{eqnarray*}
(f')^{(1)}_{i,j,k+n} &=& (h_{ij})^N (h_{jk})^N S S^{-1} (h_{ik})^{-N} c(h_{ij}, e^{i \theta})^{-1} c(h_{ik}, e^{i \theta}) \\
&=& c(h_{jk}, e^{i \theta})^{-1} c(h_{ik}, e^{i \theta}) = (\delta c^{-1})_{i,j,k+n}.
\end{eqnarray*}
The cyclic permutations are also components of the coboundary $\delta c^{-1}$. Then we have proved that $f'' = (\delta c^{-1}) (f')$. Therefore, in cohomology, the gerbe cocycle $f'$ is equal to \eqref{gerbecocycle2}. 

If we adjust the operator $S$ on the left side and use the same strategy as above, we find that the cocycle $f$ in \eqref{f} is equivalent to $f'$ in cohomology.\5 $\square$ \3 

\noindent \textbf{3.4.} We consider $\mathbb{R} \times M$ as a covering space over $\mathbb{T} \times M$ under the standard covering projection $(y,x) \mapsto ([y],x)$. Consider a lift of $\textbf{PF}$ to the covering $\mathbb{R} \times M$. The Dixmier-Douady class trivializes and therefore we can fix the phases to define a Hilbert bundle $\textbf{F}$ with a structure group $U(\mathcal{F})$ on $\mathbb{R} \times M$.  One uses the lifted transition functions to glue the fibres together. The vacuum of the above construction can be further twisted by a complex vector bundle of finite rank. Let $\xi$ denote a rank $n$ vector bundle over $M,$ extended trivially to $\mathbb{R} \times M.$ Replace next the Fock bundle $\textbf{PF}$ by $\textbf{P}(\textbf{F} \otimes \xi)= \textbf{PF}_{\xi}$ and $\textbf{F}$ by $\textbf{F} \otimes \xi = \textbf{F}_{\xi}$. Define now the action of the shift operator $S$ on the tensor product as $S\otimes 1$.

\section{Twisted K-Theory on $\mathbb{T} \times M$}

\noindent In this section we introduce a construction of twisted K-theory elements in the case of a decomposable twisting. The Fredholm operators in this model can be viewed as Dirac operators on the loop group $L \mathbb{T}$ and for this reason we begin the discussion by defining spinor modules associated to the complexified Clifford algebra of $L \mathbb{T}$. \3

\noindent \textbf{4.1.} Consider a real Clifford $*$-algebra, $\mathrm{cl}(L \mathbb{T})$ generated by $\psi_n, n \in \mathbb{Z}$ subject to the relations
\begin{eqnarray*}
[ \psi_n, \psi_m ]_+ =2 \delta_{n,-m}, \5 \psi_n^* = \psi_{-n}. 
\end{eqnarray*}
We can fix an irreducible vacuum representation of $\mathrm{cl}(L \mathbb{T})$ such that the circle group $\mathbb{T}$ acts on the vacuum
$\eta_0$ by the identity homomorphism. The operators $\psi_i$ with $i <0$  annihilate the vacuum and  the vectors $\psi_i$ with $i>0$  are used to generate the basis from the vacuum subspace. We fix the sign of $\psi_0$ such that $\psi_0 \eta_0 = \eta_0.$ We denote by $\mathcal{H}_s$ this representation.
 
On the parameter space $\mathbb{R} \times M $ we define a trivial infinite dimensional spinor bundle $\textbf{S} = \mathcal{H}_s \times \mathbb{R} \times
M$. This is pushed down to a trivial bundle over $\mathcal{H}_s \times \mathbb{T} \times M$, to be denoted by the same symbol $\textbf{S}$. Then we form a $PU$-bundle $\textbf{P}(\textbf{S} \otimes \textbf{F}_{\xi})$ over $\mathbb{T} \times M$. This does not change the topological type of the gerbe: the characteristic class of a gerbe is invariant under tensoring with trivial gerbes. We also have the Hilbert bundle $ \textbf{S} \otimes \textbf{F}_{\xi}$ over
$\mathbb{R} \times M$. 
 
We define a family of supercharge operators $Q : \mathbb{R} \times M \rightarrow  \bold{Fred}^{(1)}_{\Psi}(\textbf{S} \otimes \textbf{F}_{\xi})$ coupled
to a constant potential $y \in \mathbb{R}$ by
\begin{eqnarray*}
Q_{y} = \sum_{n \in \mathbb{Z}} \psi_n \otimes e_{-n} + y \psi_0  \otimes \textbf{1}
\end{eqnarray*}
where the operators $e_n$ define a projective unitary representation of the loop algebra $\mathfrak{lt}$ (Lie algebra of $L\mathbb{T}$) on the fibres of $\mathbf{F}$. More precisely we can write 
\begin{eqnarray*}
e_n = \sum_i : a^*(v_{n+i})a(v_i):.
\end{eqnarray*}
These operators are globally defined on $M$ (and therefore we have dropped the argument
$x\in M$ from the family $Q$). Initially we need to fix a phase from the twisting bundle $\lambda$ to make $a^*(v_n)$ and $a(v_m)$ well-defined but since the first one is linear whereas the second one is antilinear these phases cancel each other. The usual 
normal ordering  $::$ is applied to make the operators well defined on the Fock
spaces; that is, 
\begin{eqnarray*}
:a^*(v_n) a(v_m): =  \left\{ \begin{array}{rl} - a(v_m)a^*(v_n)  & \text{if $n=m < 0$}, \\
a^*(v_n) a(v_m) & \text{otherwise} 
\end{array} \right.
\end{eqnarray*}
Notice that $Q$ is unbounded and self-adjoint since $\psi^*_n = \psi_{-n}$ and $e_n^* = e_{-n}$. The finite particle states give a dense domain for $Q$. These are the states which have finite charge as well as finite fermion number: the number of operators $\psi_n$ needed to generate the state from the vacuum. \3

\noindent \textbf{Proposition.} The kernel bundles of the family $Q: \mathbb{R} \times M \rightarrow \bold{Fred}^{(1)}_{\Psi}$ localize over the submanifolds $\{k\} \times M$ for all $k \in \mathbb{Z}$. \3

\noindent Proof. The square of $Q$ is given by
\begin{eqnarray*}
Q_{y}^2 = \sum_{n > 0} n \psi_n \psi_{-n} + 2\sum_{n > 0} e_n e_{-n} + e_0^2 + 2 y e_0 + y^2 := l_0^{s} + l_0^f + (e_0 + y)^2.
\end{eqnarray*}
The operators $l_0^s$ and $l_0^f$ are positive. It is known from the representation theory of $\mathfrak{lt}$ that 
\begin{eqnarray*}
[l_0^f, S] = 0 = [l_0^s, S].
\end{eqnarray*}
Moreover, the zero modes of of $l_0^s$ and $l_0^f$ are given by $\eta_0 \otimes S^n|0 \rangle$ for all $n \in \mathbb{Z}$ and for all $1 \leq i \leq n$. The operator $e_0$ counts the charge of the state and $S e_0 S^{-1} = e_0 - 1$. It follows
\begin{eqnarray*}
Q^2_{y} (\eta_0 \otimes S^n |0 \rangle) = (n + y)^2  (\eta_0 \otimes S^n |0 \rangle)
\end{eqnarray*}
We conclude that there is a zero subspace of $Q_{y}$ if and only if $y \in \mathbb{Z}$. \5 $\square$ \3

It is also evident from the proposition that on the slice $\{k\} \times M$, with $k \in \mathbb{Z}$, there is a vacuum subbundle with one dimensional fibres and which is given locally on $\{k\} \times M$ by $\eta_0 \otimes S^k |0 \rangle$.  

The shift operators $S$ act on the supercharge by conjugation such that
\begin{eqnarray*}
S Q_{y} S^{-1} = Q_{y}.
\end{eqnarray*}
This holds since $S$ commutes with the operators $e_n$ unless $n = 0$, and as noted above, $e_0$ counts the charge. Therefore, if we set $y = \phi / 2 \pi$:
\begin{eqnarray*}
\hat{g}^{(1)}_{i,j+n} Q_{\frac{\phi}{2\pi}} (\hat{g}^{(1)}_{i,j+n})^{-1} = Q_{\frac{\phi}{2 \pi} - 1}.
\end{eqnarray*}
Consequently, we can realize $Q$ as locally defined families
on $\mathbb{T}\times M$, 
 $Q^i: U_i \rightarrow  \bold{Fred}^{(1)}_{\Psi}$, glued together by an adjoint action of a {\v C}ech-cocycle of the gerbe defined in Section 3. These are indeed Fredholm families: notice that the 1-dimensional kernel subbundles over the slices $\{k\} \times M$ are identified over $[0] \times M$ in $\mathbb{T} \times M$. The local approximated sign families
\begin{eqnarray*}
F^i := \frac{Q^i}{\sqrt{1 + (Q^i)^2}}. 
\end{eqnarray*}
define a section in the bundle of bounded self-adjoint Fredholm operators associated to the gerbe. It remains to check the continuity issues to conclude that the section $F$ defines a class in the twisted K-theory group $K^1(\mathbb{T} \times M, \sigma)$ where $\sigma$ is the decomposable class. \3

\noindent \textbf{Theorem.} \textit{The section $F$ defines an element in the group $K^1(\mathbb{T} \times M, \sigma)$.} \3

\noindent Proof. To prove the continuity it is sufficient to work on the cover $\mathbb{R} \times M.$
The  family $1 - (F_y)^2 = (1 + (Q_y)^2)^{-1}$ given by 

\begin{eqnarray*}
1 - (F_y)^2 = \frac{1}{1 + l_0^{s} + l_0^f + (e_0 + y)^2}
\end{eqnarray*}

consists of compact operators since the finite particle subspaces are finite dimensional and, as the charge in the Fock space and the fermion number in the spinor space go to infinity, the eigenvalue of $1 - (F_y)^2$ goes to $0$ for any  $y$. Moreover, $1 - (F_y)^2$ is norm continuous  since

\begin{eqnarray*}
||\frac{1}{1 + (Q_y)^2} - \frac{1}{1 + (Q_{y'})^2}||  &=&  || \frac{Q_{y'} + Q_y}{(1 + (Q_y)^2)(1 + (Q_{y'})^2)}(Q_{y'} - Q_y )|| \\
& = & ||\frac{Q_{y'}+ Q_y}{(1 + (Q_y)^2)(1 + (Q_{y'})^2)}|| \cdot || (y' - y) \psi_0 \otimes \textbf{1}||   \\
& \leq &   M ||  (y'-y)\psi_0 \otimes \textbf{1}||
\end{eqnarray*}
where $M \in \mathbb{R}$. The inequality holds since $F_y (1 + (Q_{y'})^2)^{-1}$ is a bounded operator for all $y,y'.$  It also follows that the section $F$ is its own parametrix. 

According to \cite{FHT2} it is sufficient to check the strong continuity of $F_y$ on an arbitrary vector $\Psi$ in the finite particle subspace which is dense in the Fock space. To this end we write
\begin{eqnarray*}
F_y - F_{y'}  &=& \Big(\frac{1}{\sqrt{1 + (Q_y)^2}} - \frac{1}{\sqrt{1 + (Q_{y'})^2}} \Big)Q_{y'} + \frac{Q_y - Q_{y'}}{\sqrt{1 + (Q_y)^2}} \\
&=& \Big(\frac{1}{\sqrt{1 + (Q_y)^2}} - \frac{1}{\sqrt{1 + (Q_{y'})^2}} \Big)Q_{y'} + \frac{1}{\sqrt{1+(Q_y)^2}}((y-y')\psi_0 \otimes  \textbf{1}).
\end{eqnarray*}
The second term is norm continuous and therefore strongly continuous. We have seen that $(1 + (Q_y)^2)^{-1}$ is norm continuous, and since the square root preserves norm continuity, we see that 
\begin{eqnarray*}
y \mapsto F_y\Psi 
\end{eqnarray*}
is continuous. \5 $\square$ \vspace{0,3cm}

\noindent \textbf{4.2.} If the ordinary K-theory groups $K^*(M)$  are known, one can use the Mayer-Vietoris sequence to study the $K$-theory on $\mathbb{T} \times M$ twisted by a decomposable $3$-cohomology class. \vspace{0,3cm}

\noindent \textbf{Theorem.} \textit{If $\sigma$ is a decomposable class associated to the line bundle $\lambda$, then the twisted K-theory groups $K^*(\mathbb{T} \times M, \sigma)$, for $* = 0,1$, are isomorphic to the group extension of}
\begin{eqnarray*}
	\{x \in K^*(M) | x = x \otimes \lambda \}\5 \textit{by} \5 K^{*-1}(M)/(K^{*-1}(M) \otimes (1-\lambda)).
\end{eqnarray*}

\noindent Proof. Consider a closed cover $\{\mathbb{T}_0, \mathbb{T}_1\}$ of $\mathbb{T}$ so that $\phi \in [0,\pi]$ in $\mathbb{T}_0$ and $\phi \in [\pi, 2 \pi]$ in $\mathbb{T}_1$. Then $\mathbb{T}_{01} = [-1] \cup [1]$. The gerbe corresponding to the decomposable cohomology class trivializes if the circle is cut. Therefore, the Mayer-Vietoris sequence in twisted K-theory gives:
\begin{eqnarray*}
\xymatrix{ & K^0(\mathbb{T} \times M, \sigma) \ar[r]^{\hspace{-1.2cm} c_0} &  K^0(\mathbb{T}_0 \times M) \oplus K^0(\mathbb{T}_1 \times M ) \ar[r]^{\hspace{1.2cm} a_0} &  K^0(\mathbb{T}_{01} \times M) \ar[d]_{b_0} \\ 
& K^1(\mathbb{T}_{01} \times M) \ar[u]_{b_1}   & \ar[l]_{\hspace{-1.2cm} a_1} K^1(\mathbb{T}_0 \times M) \oplus K^1(\mathbb{T}_1 \times
M )  & \ar[l]_{\hspace{1.2cm} c_1}  K^1(\mathbb{T}  \times M,\sigma) }.
\end{eqnarray*}
Thus, there are the following group isomorphism  
\begin{eqnarray*}
K^{*+1}(\mathbb{T} \times M, \sigma) &\simeq& (K^*(\mathbb{T}_{01}\times M)/\mathrm{Im}(a_*)) \oplus_{\zeta}  \mathrm{Im}(c_{*+1}) \\
&\simeq& (K^*(M)^{\oplus 2} / \mathrm{Im}(a_*)) \oplus_{\zeta} \mathrm{Ker}(a_{*+1})
\end{eqnarray*}
which is a group extension of $\mathrm{Ker}(a_{*+1})$ by $K^*(M)^{\oplus 2} / \mathrm{Im}(a_*)$ associated
to some cocycle $\zeta$ in the group cohomology.

We need to apply the tensor product operation by the bundle $\lambda$ over $M$ when we transform from $\mathbb{T}_1 \times M$ to $\mathbb{T}_0 \times M$ in $[1] \times M$. This will have effect on the glueing maps $a_*$. Consider a class $(x,y) \in K^*(\mathbb{T}_{0} \times M) \oplus
K^*(\mathbb{T}_{1} \times M)$ for $* \in \{0,1\}$. The gluing maps $a_*$ are defined by
\begin{eqnarray}\label{eq:ahom}
a_*(x,y) = (x-y, x-y\otimes \lambda)
\end{eqnarray}
where the first component on the right side is a group element in
$K^*([-1] \times M)$ and the second in $K^*([1] \times M)$. The tensor product is defined
by the usual ring structure in the ordinary $K$-theory. From \eqref{eq:ahom} we obtain
\begin{eqnarray*}
K^*(M)^{\oplus 2}/\mathrm{Im}(a_*) = K^*(M)/ (K^*(M) \otimes (1-\lambda)).
\end{eqnarray*}
Moreover,
\begin{eqnarray*}
 \mathrm{Ker}(a_{*+1}) = \{x \in K^{*+1}(M) | x = x \otimes \lambda \}
\end{eqnarray*}
This proves the theorem. \5 $\square$ \vspace{0,3cm}

In general, the above formula solves the twisted K-theory groups only up to a group extension. However, the group extension is necessarily trivial unless both groups in the extension problem have torsion subgroups. For example, when $M= S^2$ is the unit sphere and $\lambda$ is the complex line bundle equal to $k$:th tensor power of the generator,
one obtains the known result $K^0(\mathbb{T} \times M, \sigma) = \mathbb{Z}$ and $K^1(\mathbb{T} \times M, \sigma) = \mathbb{Z} \oplus \mathbb{Z}_k$ and  \cite{BCMMS}.
For $M= \mathbb{T}^2$ the corresponding groups are $\mathbb{Z}^{\oplus 2}$ and $\mathbb{Z}^{\oplus 2} \oplus \mathbb{Z}_k$.

\section{Differential twisted K-theory}
 
\textbf{5.1.} Next we study the twisted $K^1$ classes of the Fredholm families introduced above. The general principle is to map the K-theory element to some cohomology theory where it becomes easier to distinguish different group elements. The usual Chern character map to de Rham forms can be applied in twisted K-theory; however, this map loses the torsion information and therefore is not sufficient in general. 

In our case the operator family  defines a class in twisted K-theory which has components both in a torsion group and in a free group $\mathbb{Z}^p.$  In order to keep track on the torsion component we need a
refinement of the twisted K-theory to a differential twisted K-theory which depends on Dixmier Douady form,
connection and curvature of the gerbe; this data is denoted by $\check{\sigma}.$  Associated to a torsion cohomology class
there is an eta form which is closed in a twisted cohomology theory $H^3(X,
H)$ \cite{CMW}.  The coboundary operator is given by $\delta = d - H$ where the three form
$H$ is the de Rham representative of the Dixmier-Douady class of the gerbe. In the case of torsion twisted K-theory, the
eta form is related to a differential twisted Chern character form by the formula
\begin{eqnarray*}
(d-H) \eta = - \ch_{\check{\sigma}}(Q) 
\end{eqnarray*}
In the following we construct the Chern character form
$\ch_{\check{\sigma}}(Q)$. \3

\noindent \textbf{5.2.} We cut the circle $\mathbb{T}$ as in the Mayer-Vietoris sequence in the proof of Theorem 4.2: $\mathbb{T}_0$ and $\mathbb{T}_1$ are the closed segments where $\phi \in [0,\pi]$ and $\phi \in [\pi, 2 \pi]$ respectively. Then we have the following maps in the Mayer-Vietoris sequence 
\begin{eqnarray*}
K^0(\mathbb{T}_0 \times M) \oplus K^0(\mathbb{T}_1 \times M) &\longrightarrow& K^0(\mathbb{T}_{01} \times M) \\ &\stackrel{\simeq}{\longrightarrow}&
K^0([-1] \times M) \oplus K^0([1] \times M).
\end{eqnarray*}
The composition sends $(x,y)$ to $(x-y) \oplus (x-y \otimes \lambda)$. Therefore, the image
of $(x \otimes \lambda, x)$ is zero in the component $[1] \times M$. We set an inclusion $i: [1] \times M \rightarrow \mathbb{T} \times
M$.  Then there is an equivalence in differential twisted K-theory
\begin{eqnarray*}
i_!^K (\xi) = i_!^K (\xi \otimes \lambda),
\end{eqnarray*}
if $\xi$ is a vector bundle and $i_!^K$ is a push forward map in the differential twisted K-theory 
\begin{eqnarray*}
i_!^K: K^0(M) \rightarrow K^1(\mathbb{T} \times M, \sigma).
\end{eqnarray*}

\noindent \textbf{5.3.} Following the general theory, \cite{CMW}, we introduce a differential twisted Chern character for the supercharge element $Q$. Similar calculations were done in the case of twisted K-theory on $SU(n)$ in \cite{CMW}. Recall that the zeros of the operators $Q_{y}$ are localized at the integer values of $y \in \mathbb{R}$ and the kernel bundle is topologically equivalent to $\xi \otimes \lambda^{ \otimes n}$ over $\{ n\} \times M $. Down on $\mathbb{T}$, the bundles are defined modulo tensoring with powers of $\lambda$, related to the choice of a local section $\mathbb{T} \to \mathbb{R}$.  However, this does not cause problems since the twisted $K^1$-theory is independent of tensor powers by $\lambda$ by 5.2. Then we have a pushforward of the class of the kernel bundle: 
\begin{eqnarray*}
	i_!^K: K^0([1] \times M) \rightarrow K^1(\mathbb{T} \times M, \sigma).
\end{eqnarray*}
Fow what follows we need to have an inclusion map so that the difference of the dimension of its range and the dimension of its domain $[1] \times M$ is even. Therefore we replace $\mathbb{T} \times M$ with $\mathbb{T} \times M \times \mathbb{R}$. Then the kernel bundle of $Q$ can be pushed by the map  
\begin{eqnarray*}
	i_!^K: K^0([1] \times M) \rightarrow K^1(\mathbb{T} \times M \times \mathbb{R}, \sigma).
\end{eqnarray*}
The K-theory group in the target is twisted by a gerbe which is extended trivially on $\mathbb{R}$.

Following the discusssion in \cite{CMW} (for twisted K-theory on $SU(n)$)  we apply the Riemann-Roch theorem in
twisted differential K-theory. In general, if $f: X \rightarrow Y$ is a closed embedding of manifolds such that $\mathrm{dim}(Y) -\mathrm{dim}(X)$ is even and the normal bundle of the embedding, $N_f \simeq f^*(TY)/TX$, has a $\mathrm{spin}^c$-structure $c_1$ and $a \in K^0(X, f^*(\sigma))$, then
\begin{eqnarray*}
	\mathrm{ch}_{\check{\sigma}}(f^K_{!}(a)) \hat{A}(Y) = f^H_* (\mathrm{ch}_{f^* \check{\sigma}}(a) e^{\frac{c_1}{2}} \hat{A}(X))
\end{eqnarray*}
where $f^H_*$ is a Gysin homomorphism in twisted cohomology theory. 

The inclusion $i: [1] \times M \rightarrow \mathbb{T} \times M \times \mathbb{R}$ is a closed embedding of manifolds with a trivial normal bundle. The cohomology class  $i^* \sigma$ and consequently $i^* \check{\sigma}$ are trivial and therefore the map $i^H_*$ is the Gysin homomorphism in the de Rham cohomology and $\mathrm{ch}_{f^* \check{\sigma}}$ is the usual Chern character. The $\hat{A}$-genus satisfies 
\begin{eqnarray*}
	\hat{A}(\mathbb{T} \times M \times \mathbb{R}) = \hat{A}(M)  \5 \hbox{and} \5 \hat{A}([1] \times M) = \hat{A}(M).
\end{eqnarray*}
Therefore, 
\begin{eqnarray*}
\ch_{\check{\sigma}}(i_!^K(\xi)) = i_*^H(\ch(\xi) \hat{A}(M)) (\hat{A}(M))^{-1} = i_*^H(\ch(\xi)).
\end{eqnarray*}
Then, by definition:
\begin{eqnarray*}
(d - H) \eta = - i_*^H(\ch(\xi)) = - i_*^H(\ch(\xi \otimes \lambda^{\otimes n})),
\end{eqnarray*}
for all $n \in \mathbb{Z}$. The pairing of $(d - H) \eta$ with a de Rham cycle is a twisted K-theory invariant. We can use this to study the dependence of $Q$ on the vector bundle $\xi$ on $M$.\3
 
\noindent \textbf{5.4.} We consider the case of a torsion class in $K^1(\mathbb{T} \times M, \sigma)$ such that the curvature of the twisting line bundle $\lambda$ is an element $k F_b$ in $H^2(M, \mathbb{Z})$ and the vacuum line bundle has the curvature $n F_b$ and $F_b$ generates a subgroup isomorphic to $\mathbb{Z}$ in $H^2(M, \mathbb{Z})$. Furthermore, we assume that the form $F_b$ localizes on a two dimensional closed oriented surface $S$ on $M$. Then, by integration we find that 
\begin{eqnarray*}
-\int_{\mathbb{T} \times S} (d - H) \eta &=&  \int_{\mathbb{T} \times S} i_*^H(\ch(\xi)) \5 \mathrm{mod} \5  n \\ 
&=& \int_{S} \ch(\xi) \5 \mathrm{mod}\5 n  \\ &=& k  \5 \mathrm{mod}\5 n.
\end{eqnarray*}
There is also a degree 1 nontorsion component in the Chern character, which in the case of a complex line bundle $\xi$ is equal to the generator in $H^1(\mathbb{T}, \mathbb{Z})$, see the computation in the next Section. On the degree $1$ level the twisted cohomology is the same as the usual de Rham cohomology since $H\wedge d\phi =0$
on $\mathbb{T} \times M.$

\section{Superconnection Analysis}  

\noindent \textbf{6.1.} We start by studying a superconnection associated to the family $Q$ over the covering space $\mathbb{R} \times M$. The connection $\nabla=\nabla_{M} \otimes 1 + 1 \otimes \nabla_{\xi}$ consists of of a connection $\nabla_{\xi}$ of the bundle $\xi$ over $M$ and a connection $\nabla_{M}$ of the twisting line bundle $\lambda$ over $M.$ The action of the connection $\nabla_{M}$  on the fermion number $n$ sector is the $n$'th tensor power of the connection in the line bundle $\lambda$; in particular, on the vacuum sector the only nontrivial piece is $\nabla_{\xi}.$
 Let us define
\begin{eqnarray*}
	\mathbb{A}_t = \sqrt{t} \chi Q + \nabla. 
\end{eqnarray*}
We have also introduced a real scaling parameter $t > 0$. We write locally  $\nabla = d +\omega$  where $\omega$ is the matrix valued connection form acting on the sections of the Fock bundle and $\hat F=\nabla^2$ is the curvature two form, composed of $\beta_M  = \nabla_{M}^2= e_0 \beta_M$ and $F_{\xi}= \nabla_{\xi}^2$. The formal symbol $\chi$ with $\chi^2=1$ is introduced since the Clifford algebra of the loop group on the circle is odd (the circle is odd dimensional).
The symbol $\chi$ is defined to commute with $Q$ and anticommute with odd differential forms. Note that the Bismut superconnection \cite{Bi} for families of Dirac
operators contains a term proportional to the curvature with a factor $1/\sqrt{t}.$ The motivation for that term is that in the limit $t\to 0$ one obtains 
from the character formula below the local Atiyah-Singer index formula. However, here we shall study the limit $t\to \infty$ and we drop this term.
We have
\begin{eqnarray*}
	d = d_{y} + d_M, \5  \hat F = d \omega + \omega^2  = e_0 \beta_M \otimes 1 +  1 \otimes F_{\xi}
\end{eqnarray*}
(we denote $y=\phi/2\pi$). The square of the superconnection is 
\begin{eqnarray*}
	\mathbb{A}_{t,y}^2 &=& t Q_y^2 + \sqrt{t} \chi( - d Q_y + [Q_y,\omega]) + \hat{F} \\
	&=& t Q_y^2 - \sqrt{t} \chi \psi_0 dy + \hat{F}. 
\end{eqnarray*}
 Using $S e_0 S^{-1} = (e_0 - 1)$ and $S Q^2_{y} S^{-1} = Q^2_{y-1} $ one gets
\begin{eqnarray*}
	\mathbb{A}_{t,y+1}^2 = S^{-1} \mathbb{A}^2_{t,y} S +  \beta_M.
\end{eqnarray*}
If we define the character of the superconnection form by 
\begin{eqnarray*}
	\Theta_{y} = \mathrm{sTr} (e^{- \mathbb{A}^2_{t,y}}), 
\end{eqnarray*}
where the supertrace $\sTr$ picks up the terms linear in $\chi$, then 
\begin{eqnarray*}
	\Theta_{y+1} = \mathrm{sTr}(e^{-\mathbb{A}^2_{t,y} -\beta_M})  = \Theta_{y} \wedge e^{-\beta_M}.
\end{eqnarray*}

The character needs to be pulled to the base $\mathbb{T} \times M$ where the twited K-theory is realized but this is not straightforward because of the nonperiodicity of the form $\Theta$ under translations in the direction of $\mathbb{R}$. One could define a new form  
\begin{eqnarray*}
	\tilde{\Theta}_{y} = e^{y\beta_M} \wedge \Theta_{y}, 
\end{eqnarray*}
so that $\tilde{\Theta}_{y+1} = \tilde{\Theta}_{y}$, pull the form $\tilde{\Theta}$ to the base and then associate a twisted cohomology class with it. The differential in the twisted cohomology is given by $d - H$ for $H = d\phi \wedge \beta_M$. We shall introduce another approach below. \3

\noindent \textbf{6.2.} We first compute the $t \rightarrow \infty$ limit of the superconnection over $\mathbb{R} \times M$, then proceed by pulling it to the base. In this limit, the character localizes over the zeros of the family $Q$. For this reason we introduce the projection operator $P$ onto the subbundle 
\begin{eqnarray}\label{P}
\eta_0 \otimes \bigoplus_{k \in \mathbb{Z}} \hat{S}^k |0 \rangle \otimes \xi. 
\end{eqnarray}
Moreover, we denote by $\delta$ the Dirac delta distribution. \3

\noindent \textbf{Proposition.} The $t \rightarrow \infty$ limit of the form $\Theta$ exists as a distributional valued form: 
\begin{eqnarray*}
 \sqrt{\pi} \delta(P e_0 P + y)dy \wedge \mathrm{tr}_{\xi} (e^{ -\hat{F} }). 
\end{eqnarray*}

\noindent Proof. Put $\mathbb{A}^2_{t} = t(Q^2 + K_t)$. We use the Volterra series 
\begin{eqnarray*}
	\Theta = \mathrm{sTr} \Big(e^{-tQ^2} + \sum_{n \geq 1} (-t)^n \int_{\triangle_n} e^{-ts_1 Q^2} K_t e^{-ts_2 Q^2} \cdots e^{-ts_n Q^2} K_t e^{-ts_{n+1} Q^2}   ds_1 ds_2 \cdots ds_{n+1}\Big). 
\end{eqnarray*}
where $\triangle_n$ denotes the standard $n$-simplex. The operator family $K_t = \frac{1}{t} \hat F - \frac{\chi}{\sqrt{t}} \psi_0 dy$ commute with $Q^2$ and we can simplify
\begin{eqnarray*}
	e^{-ts_1 Q^2} K_t e^{-ts_2 Q^2} \cdots e^{-ts_n Q^2} K_t e^{-ts_{n+1} Q^2} = K_t^n e^{-t Q^2}
\end{eqnarray*}
The forms $ \chi dy$ and $\hat F$ commute. Thus, 
\begin{eqnarray*}
	 K_t^n = \sum_{k = 0}^{n} \binom{n}{k} \Big(\frac{-\chi \psi_0 dy}{\sqrt{t}}  \Big)^{n-k} \wedge \Big(\frac{\hat{F}}{t}\Big)^{k} = - n \frac{\chi \psi_0 dy}{\sqrt{t}} \wedge \Big( \frac{\hat{F}}{t}\Big)^{n-1} + \Big(\frac{\hat{F}}{t}\Big)^{n},
\end{eqnarray*}
 A volume of the simplex is $1/n!$ and therefore we get 
\begin{eqnarray*}
	\Theta &=& \mathrm{sTr} \Big(e^{-tQ^2} +  \sqrt{t} \sum_{n \geq 1} \frac{\chi \psi_0 dy \wedge (-\hat{F})^{n-1}}{(n-1)!} e^{-tQ^2} + \sum_{n \geq 1} \frac{(-\hat{F})^n}{n!} e^{-tQ^2} \Big)  \\ 
	&=& \mathrm{Tr} \Big( \sqrt{t} \sum_{n \geq 1} \frac{\psi_0 dy \wedge (-\hat{F})^{n-1}}{(n-1)!} e^{-tQ^2} \Big). 
\end{eqnarray*}
Recall that
\begin{eqnarray*}
	Q^2 = l_0^s + l_0^f + (e_0 + y)^2. 
\end{eqnarray*}
We use the asymptotic expansion for the positive operator $e^{-t(l_0^s + l_0^f)}$ as $t \rightarrow \infty$. In this limit, the operator $e^{-tQ^2}$ converges to zero outside the subbundle \eqref{P}. The following formulas hold for the Dirac measure
\begin{eqnarray*}
	\delta(\phi - a) = \lim_{t \rightarrow \infty} \sqrt{\frac{t}{\pi}} e^{-t(\phi - a)^2}, \5 \lim_{t \rightarrow \infty} \frac{1}{t^p}\sqrt{\frac{t}{\pi}} e^{-t(\phi - a)^2} = 0\5 p \in \mathbb{N}.
\end{eqnarray*}
Therefore 
\begin{eqnarray*}
	\lim_{t \rightarrow \infty} \Theta &=&  \sqrt{\pi} \mathrm{Tr} \Big(P \psi_0  \delta( e_0 + y) dy \wedge e^{ -\hat{F} } P \Big) \\ &=& \sqrt{\pi} \delta(P e_0 P + y) dy \wedge \mathrm{tr}_{\xi} (e^{ -\hat{F} }). \5 \square
\end{eqnarray*}

Choose a section $\psi$ for the standard projection $\mathbb{R} \times M \rightarrow \mathbb{T} \times M$. The section $\psi$ is only locally defined and, as discussed above, the local differential forms do not agree on all the intersections because of the nonperiodicity of $\Theta$ under translations. The nonperiodicity can be fixed as follows. Define the following normal subgroup in the de Rham cohomology group $H^*(\mathbb{T} \times M, \mathbb{Q})$: 
\begin{eqnarray*}
N =  \varphi \Big[\frac{d \phi}{2 \pi}\Big] \wedge \ch(K^0(M)) \wedge \Big(1 - \text{ch}(\lambda) \Big) 
\end{eqnarray*}
where $[d \phi]$ is the class of the form $d \phi$ and ch denotes the Chern character $K^0(M) \rightarrow H^{\text{even}}(M)$. The symbol $\varphi$ denotes the standard normalization: 
\begin{eqnarray*}
   \varphi : \Lambda_{\mathbb{C}}^{2k + 1}(M) \rightarrow \Lambda^{2k + 1}_{\mathbb{C}}(M), \5 \varphi(\Omega) = (2 \pi i )^{-k} \Omega.
\end{eqnarray*}
for all $k \in \mathbb{N}_0$. The normalization is applied to make the character get values in the rational cohomology of $\mathbb{T} \times M$. \3

\noindent \textbf{Theorem.} \textit{For all $t > 0$, the pullback form $\frac{1}{\sqrt{\pi}} \varphi(\psi^*(\Theta_t))$ defines a class in the quotient group}
\begin{eqnarray}\label{character}
 \varphi \Big[ \frac{d \phi}{2 \pi} \wedge \ch(\xi) \Big] \in H^{\text{odd}}(\mathbb{T} \times M, \mathbb{Q}) / N
\end{eqnarray}
\textit{The class is independent on the choice of $\psi$ and $t$}. \3

\noindent Proof. The $t$ independence of the character follows from the standard transgression formula, \cite{BGV04}.

Consider the distribution valued form $\delta(P e_0 P + y ) d y$ in $\mathbb{R}$ of the proposition. The pullback of this form to $\mathbb{T}$ by a section of the projection $\mathbb{R} \rightarrow \mathbb{T}$ has its support at $[0] \in \mathbb{T}$ and the pairing with the fundamenal cycle over $\mathbb{T}$ gives $1$. Therefore, by continuity, the distributional form over $\mathbb{T}$ is a limit of a differential form on $\mathbb{T}$ which is represented in the de Rham cohomology by $\frac{d \phi}{2 \pi}$. Then the pullback form $\frac{1}{\sqrt{\pi}} \varphi(\psi^*(\Theta_t))$ is represented in cohomology by the form \eqref{character}: the nonperiodicity is cancelled since in the quotient group the form is defined up to multiplets of $\ch(\lambda) = e^{-\beta_M}$. For this same reason, the choice of the section $\psi$ does not have effect on the character in $H^{\text{odd}}(\mathbb{T} \times M, \mathbb{Q}) / N$. \5 $\square$  \3

The analysis above proves that the twisted $K^1$ class associated to the family $Q$ and the vacuum vector bundle $\xi$ are distinguished by the Chern character of $\xi$ evaluated in the quotient. In the case of torsion twisted K-theory class this result is compatible with the analysis in 5.4. \3

\noindent \textbf{6.3.} Consider the case $\mathrm{dim}(M) = 2$, then the character \eqref{character} is represented in the quotient by  
\begin{eqnarray*}
	 \Big[ \frac{d \phi}{2 \pi} \wedge (\mathrm{rk}(\xi) - \mathrm{tr}_{\xi}  \frac{F_{\xi}}{2 \pi i}) \Big] \5 \text{mod} \5 \frac{d \phi}{2 \pi}  \wedge \frac{\beta_M}{2 \pi i}
\end{eqnarray*}
where $\text{tr}_{\xi} (F_{\xi})$ in the case dim $M=2$ is an integer $n$ times the curvature $F_b$ of the basic line bundle over $M$. In the case $M= S^2$ and $\beta_M = k F_b$ and $F_{\xi} = n F_b$, the operator family defines the twisted $K^1$-group element 
\begin{eqnarray*}
([n], \text{rank}(\xi)) \in K^1(\mathbb{T} \times S^2, \sigma) \simeq \mathbb{Z}_k \oplus \mathbb{Z}
\end{eqnarray*}
where the class $\sigma$ is the decomposable class associated to $\beta_M$.

\bibliographystyle{plain}

\end{document}